\documentclass[12pt]{amsart}
\usepackage[all]{xy}
\usepackage{amssymb}

\newcommand{\bbf}{\mathbb{F}}
\newcommand{\bbn}{\mathbb{N}}
\newcommand{\bbp}{\mathbb{P}}
\newcommand{\bbq}{\mathbb{Q}}
\newcommand{\bbr}{\mathbb{R}}
\newcommand{\bbz}{\mathbb{Z}}
\newcommand{\Hom}{\mathrm{Hom}}
\newcommand{\Sym}{\mathrm{Sym}}
\newcommand{\Gal}{\mathrm{Gal}\,}

\newcommand{\cha}{\mathrm{char}}

\newtheorem{thm}{Theorem}

\newtheorem{lem}[thm]{Lemma}
\newtheorem{prop}[thm]{Proposition}

\newtheorem{defnn}[thm]{Definition}

\newtheorem{remarkk}[thm]{Remark}

\newtheorem{examplee}[thm]{Example}

\title{A Weak Chevalley-Warning Theorem for Quasi-finite Fields}
\author{Bo-Hae Im and Michael Larsen}
\date{February 21, 2008}
\address{Department of Mathematics, Chung-Ang University, 221, Heukseok-dong, Dongjak-gu, Seoul, 155-756, South Korea}\email{bohaeim@gmail.com}
\address{Department of Mathematics, Indiana University, Bloomington,
Indiana 47405, USA} \email{mjlarsen@indiana.edu}
\thanks{Michael Larsen was partially supported by NSF grant DMS-0354772.}

\begin{document}
\begin{abstract} There exists a function $f\colon \bbn\to\bbn$ such that for every positive integer $d$,
every quasi-finite field $K$ and every projective hypersurface $X$ of degree $d$ and dimension
$\ge f(d)$, the set $X(K)$ is non-empty.  This is a special case of a more general result about
intersections of hypersurfaces of fixed degree in projective spaces of sufficiently high dimension over fields with finitely generated Galois groups.

\end{abstract}
\maketitle

\section{Introduction}

The Chevalley-Warning theorem \cite[I,~Th.~3]{Se2} asserts that every finite field $K$ is $C_1$, i.e., for every $n\in\bbn$ and every degree $n$ hypersurface $X$ of dimension $\ge n$, the set $X(K)$ is non-empty.  One might ask whether the Chevalley-Warning theorem extends to all quasi-finite fields (i.e., perfect fields with Galois group isomorphic to $\hat\bbz$).   Ax \cite{Ax}, answering a question of Serre \cite[II,~\S3]{Se1},  gave an example of a quasi-finite field which is not $C_n$ for any $n$.   

In this note, we prove a weak version of Chevalley-Warning, as follows:

\begin{thm}
\label{qf}
There exists a function $f\colon \bbn\to\bbn$ such that for every positive integer $d$,
every quasi-finite field $K$, and every projective hypersurface $X$ of degree $d$ and dimension
$\ge f(d)$, the set $X(K)$ is non-empty.
\end{thm}

The proof can in principle be used to produce a function $f$, but it has very rapid growth. 
Of course, we have already seen that $f(d)$ cannot have polynomial growth.
 
Theorem~\ref{qf} is a special case of the following more general result.
\goodbreak

\begin{thm}
\label{Linear}
There exists a function $h\colon \bbn^4\to\bbn$ such that given:
\begin{itemize}
\item any perfect field $K$, not formally real, such that $G_K$ has a generating set with $g$ elements,
\item any positive integers $d$, $e$, and $k$,
\item any sequence $d_1,\ldots,d_e$ of positive integers $\le d$,
\item any vector space $V$ over $K$ of dimension $n+1 > h(d,e,g,k)$,
\item any sequence of forms $F_1\in \Sym^{d_1}V^*,\ldots, F_e\in\Sym^{d_e}V^*$,
\end{itemize}
there exists a subspace $W\subset V$ of dimension $k$ on which all the forms $F_i$
are identically zero.
\end{thm}

Specializing to $k=1$, we get the following:

\begin{thm}
\label{General}
There exists a function $h\colon \bbn^3\to\bbn$ such that if $K$ is any perfect field, not formally real, such that $G_K$ has a generating set with $g$ elements; $d$ and $e$ are positive integers; and $X$ is an intersection of $e$
degree $\le d$ hypersurfaces in $\bbp^n$ for $n\ge h(d,e,g)$, then $X(K)$ is non-empty.
\end{thm}

We remark that there are interesting examples of fields $K$ for which $G_K$ is finitely generated but not
abelian.  For instance, it is known \cite[Th.~7.5.10]{NSW}
that every $p$-adic field has a finitely generated Galois group.
For $K=\bbq_p$, $p>2$, we can take $g=4$.  For this family of fields, 
Theorem~\ref{General} is due to Schmidt  \cite{Sc}.

Specializing Theorem~\ref{General} to the case $e=g=1$, we obtain Theorem~\ref{qf}, since by
definition, quasi-finite fields are perfect, and a quasi-finite field cannot be formally real.
Indeed, the Brauer group of a quasi-finite field is trivial \cite[XIII, Prop.~5]{Se3}.  Therefore,
the Severi-Brauer curve $x^2+y^2+z^2=0$ (\cite[X, \S7, Ex. (e)]{Se3}) has a rational point over $K$,
which means that $-1$ is a sum of two squares in $K$.  Thus $K$ satisfies the hypotheses of 
Theorem~\ref{General}.

The idea of the proof of Theorem~\ref{Linear} is to use a polarization argument due to Brauer \cite{Br} to reduce to the case of diagonal forms.  By a Galois cohomology argument, we can further reduce to the case of Fermat hypersurfaces.  These can be treated using identities introduced by Hilbert in his work on Waring's problem.  

We would like to thank Kiseop Park and Olivier Wittenberg for calling our attention to relevant literature.

We begin with Fermat hypersurfaces

$$F_{d,n}:\ x_0^d+x_1^d+\cdots+x_n^d = 0.$$

\begin{thm}
If $d$ is a positive integer and $K$ is a field, not formally real,
such that $K_d := K^\times/(K^\times)^d$ is finite, then $F_{d,n}(K)$ is non-empty
whenever $n\ge|K_d|$.
\end{thm}

\begin{proof}
Let
$$\Sigma_{d,r} := (\underbrace{(K^\times)^d + \cdots + (K^\times)^d}_r)\cup\{0\}.$$
Then $\Sigma_{d,r}$ is stable by multiplication by $(K^\times)^d$ and is therefore characterized by 
the image in $K_d$ of its non-zero elements.  As $\Sigma_{d,r}\subset \Sigma_{d,r+1}$, it follows that the
sequence must stabilize, and every sum of $d$th powers of elements of $K$ can be expressed as a sum of at most $|K_d|$ such powers.  

It remains to show that for every positive integer $d$, $-1$ is a sum of $d$th powers.
If the characteristic of $K$ is positive, this is obvious, so we assume that $\cha\ K = 0$.
By hypothesis, $-1$ is a sum of squares in $K$.  This implies that $\Sigma_{2,r} = K$ for $r\gg 0$.  Indeed,
the condition $a_1^2+\cdots+a_n^2 = -1$ implies 
$$\Bigl(\frac{c+1}2\Bigr)^2 + \sum_{i=1}^n \Bigl(a_i\Bigl(\frac{c-1}2\Bigr)\Bigr)^2 = c.$$
For odd $d$, it is obvious that $-1$ is a sum of $d$th powers.  If for some $d$ and some $r$,
$-1\in \Sigma_{d,r}$, then we can write
$$-1 = \sum_{i=1}^r\Bigl(\sum_{j=1}^{|K_2|} a_{i,j}^2\Bigr)^d.$$
By Hilbert's identity \cite[Th.~3.4]{Na}, this implies that $-1$ is a sum of $2d$th powers
of certain $\bbq$-linear combinations of the $a_{i,j}$.
The theorem follows by induction on the largest power of $2$ dividing $d$.
 
\end{proof}

Next we consider diagonal hypersurfaces in general.

\begin{lem}
\label{minimum-generators}
Let $G$ be a profinite group and $H < G$ an open subgroup.  If $G$ can be topologically generated by $g$ elements, then $H$ can be topologically generated by 
$1+[G:H](g-1)$ elements.
\end{lem}

\begin{proof}
By compactness, it suffices to prove this when $G$ and $H$ are finite groups.
Let $G'$ denote a free group on $g$ elements and $\pi\colon G'\to G$ a surjective homomorphism.
Let $H' := \pi^{-1}(H)$.  Thus $H'$ is a subgroup of $G'$ of index $[G:H]$.
Identifying $G'$ with the fundamental group of a join $X$ of $g$ circles, the quotient of the universal
cover of $X$ by $H'$ is a covering space of $X$ of degree $[G':H']$ and therefore a finite, connected, 1-dimensional CW-complex $Y$.  Thus $H'$ is free, and the number of its generators is the rank of
$H_1(Y,\bbz)$.  We have
$$\chi(Y) = [G':H']\chi(X) = [G:H]\chi(X) = [G:H](g-1),$$
where $\chi$ denotes the Euler characteristic.  Thus, $H'$ has $1+[G:H](g-1)$ generators, and the same is true of its quotient $H$.
\end{proof}

\begin{prop}
If $K$ is a perfect field such that $G_K$ can be generated by $g$ elements and $d$ is a positive integer,  then
$$|K_d|\le d^{dg+1}.$$
\end{prop}

\begin{proof}
Let $K_d$ denote the extension of $K$ generated by $\mu_d$, the group of solutions of $x^d = 1$ in 
$\bar K$.
We have Kummer isomorphisms 
$$K^\times/(K^\times)^d\cong H^1(K,\mu_d)$$
and
$$K_d^\times/(K_d^\times)^d\cong H^1(K_d,\mu_d)\cong \Hom(G_{K_d},\mu_d)$$

Now, $G_{K_d}$ is the kernel of the homomorphism $G_K\to \Gal(K_d/K)$ whose image has
order $\le d-1$.  
Applying Lemma~\ref{minimum-generators}
to $G_{K_d}\subset G_K$,  we see that 
$$|\Hom(G_{K_d},\mu_d)|\le d^{(d-1)(g-1)+1}.$$
The order of the cohomology group $H^1(\Gal(K_d/K),\mu_d)$ is bounded above by the number of $1$-cochains.
By the inflation-restriction sequence, 
$$|K^\times/(K^\times)^d| \le |H^1(\Gal(K_d/K),\mu_d)|\cdot |H^1(K_d,\mu_d)|\le d^{dg+1}.$$
The proposition follows.
\end{proof}

We can now prove Theorem~\ref{Linear}.

\begin{proof}
By \cite[Th.~C]{Br}, it suffices to find a bound $N$, depending only on $d$ and $g$, such that every
diagonal homogeneous form of degree $d$ in $K$ in more than $N$ variables has a solution in $K$. 
Let $n = d^{dg+1}$, $N = n^2$, and 
$$F:=a_0 x_0^d+\cdots+a_{N} x_{N}^d$$ 
a diagonal form.  If any coefficients $a_i$ is zero, then $F$ has an obvious non-trivial solution.
Otherwise at least  $n+1$ of the $a_i$ must belong to a single
coset $a (K^\times)^d$.  Assuming without loss of generality that $a_0,\ldots,a_n\in a(K^\times)^d$ and letting $(b_0,\ldots,b_n)$ denote a non-trivial solution of $x_0^d+\cdots+x_n^d=0$ in $K$,
the non-zero vector
$$(b_0,\ldots,b_n,0,\ldots,0)$$
is a solution of $F$.
\end{proof}

We conclude with two examples, to show that the hypotheses on $K$ are really needed.

If $K = \bbr$ and $F:= x_1^2+x_2^2+\cdots+x_n^2$, there are no non-trivial solutions regardless of how large $n$ may be.

If $K$ is a separable closure of the infinite-dimensional function field $\bbf_p(t_1,t_2,\ldots)$, then
$$F := t_1 x_1^p + t_2 x_2^p + \cdots + t_n x_n^p = 0$$
has no solutions regardless of how large $n$ may be since the $t_i$ are linearly independent over
$\bbf_p(t_1^p,t_2^p,\ldots)$.

\end{document}